\documentclass[12pt]{article}%31.01.2016
\usepackage{amsmath}
\usepackage{amssymb}
\usepackage{dsfont}
\usepackage{cite}
\newsymbol \blackbox 1004
\newcommand{\eh}{\hfill}\newlength{\sperr}

\newenvironment{proof}{{\settowidth{\sperr}{\bf\rm
Proof}%
\par\addvspace{0.3cm}\noindent\parbox[t]{1.3\sperr}
{\bf\rm P\eh r\eh o\eh o\eh f\eh }%
}}{\nopagebreak\mbox{}
$\blackbox$\par\addvspace{0.3cm}}

\def\nn{\nonumber}
\def\a{\alpha}

\def\g{\gamma}

\def\vk{\varkappa}

\def\Lam{\Lambda}
\def\s{\sigma}

\def\t{\theta}

\def\vp{\varphi}
\def\vt{\vartheta}
\def\ve{\varepsilon}
\def\wh{\widehat}
\def\wt{\widetilde}

\def\BC{{\mathbb C}}
\def\BR{{\mathbb R}}
\def\BN{{\mathbb N}}

\def\cla{{\mathcal A}}
\def\clb{{\mathcal B}}
\def\clc{{\mathcal C}}

\def\cls{\mathcal{S}}
\def\clt{\mathcal{T}}

\def\im{{\rm Im\ }}
\def\spa{{\rm Span}}
\def\diag{\mathrm{diag}}
\newcommand{\E}{\mathrm{e}}
\newcommand{\I}{\mathrm{i}}

\newtheorem{Pa}{Paper}[section]
\newtheorem{Tm}[Pa]{{\bf Theorem}}
\newtheorem{La}[Pa]{{\bf Lemma}}
\newtheorem{Cy}[Pa]{{\bf Corollary}}
\newtheorem{Rk}[Pa]{{\bf Remark}}

\newtheorem{Dn}[Pa]{{\bf Definition}}
\newtheorem{Nn}[Pa]{{\bf Notation}}
\newtheorem{Pn}[Pa]{{\bf Proposition}}

\title{Skew-selfadjoint Dirac systems: stability of the procedure of explicit solving the inverse problem}

\author{B. Fritzsche, B. Kirstein, I.Ya. Roitberg and A.L. Sakhnovich}

\date{}
\begin{document}
\maketitle

\begin{abstract} Procedures to recover explicitly discrete and continuous skew-selfadjoint Dirac systems on semi-axis
 from rational Weyl matrix functions are considered. Their stability is shown. Some new facts on asymptotics of pseudo-exponential potentials
 (i.e., of explicit solutions of inverse problems)
 are proved as well.
 GBDT version of B\"acklund-Darboux transformation, methods from system theory
and results on algebraic Riccati equations are used for this purpose.
\end{abstract}



\section{Introduction}
\setcounter{equation}{0}
Skew-selfadjoint Dirac system on the semi-axis has the form
\begin{align} &       \label{1.1}
y^{\prime}(x, z )= (\I z j+jV(x))y(x,
z ), \qquad
x \geq 0 \quad (z\in \BC),
\end{align} 
where
\begin{align} &   \label{1.2}
j = \left[
\begin{array}{cc}
I_{m_1} & 0 \\ 0 & -I_{m_2}
\end{array}
\right], \hspace{1em} V= \left[\begin{array}{cc}
0&v\\v^{*}&0\end{array}\right],  \quad m_1+m_2=:m,
 \end{align} 
$y^{\prime}=\frac{d}{dx}y$, $\, I_{m_k}$ is the $m_k \times m_k$ identity
matrix, $v(x)$ is an $m_1 \times m_2$ matrix function and $\BC$ stands for the complex plane.
In this paper, we assume that the {\it potential} $v$ is bounded, that is, 
\begin{align} &   \label{1.2!}
\sup_{x\in [0,\infty)}\|v(x)\| \leq M
 \end{align} 
for some $M>0$. Here $\|\cdot \|$  is the $l^2$-induced matrix norm.

Discrete skew-selfadjoint Dirac system  is given (see \cite{KS, FKKS}) by the formula:
\begin{equation} \nonumber
y_{k+1}(z)=\left(I_m+ \frac{\I}{ z}
C_k\right)
y_k(z),  \quad C_k=U_k^*jU_k \quad \left( k \in \BN_0 \right),
\end{equation}
where the matrices $U_k$ are unitary, $j$ is defined in \eqref{1.2}, and  $\BN_0$ stands for the set of non-negative integers.

Inverse spectral problems to recover systems from spectrum or from  Weyl functions are usually nonlinear and unstable and
the cases of stability of the procedure are especially  interesting. Here, we deal with the inverse problem to recover systems from  Weyl functions.
A procedure  of explicit solving the inverse problem for continuous
selfadjoint Dirac system was worked out in  \cite{FKRS2013, GKS6}, and the stability of this procedure was recently studied in \cite{ALSstab}. Skew-selfadjoint Dirac systems are in many respects as important as selfadjoint ones but they
present also some additional difficulties being non-selfadjoint. The discrete case is in many respects even more complicated
than the continuous. General-type inverse problems to recover skew-selfadjoint Dirac systems from Weyl functions
were studied in \cite{CG2, FKRS2012, SaA021, SaA14, SaSaR}.  Explicit solutions of inverse problems are often obtained in
a different (from general-type solutions) way, using Crum-Krein methods \cite{Cr, KrChr}, commutation methods \cite{D, Ge, GeT, T00} and
various versions of B\"acklund-Darboux transformations (see, e.g., \cite{Ci, MS, SaA1, SaSaR, ZM} and numerous references therein).
We consider here the GBDT (generalized B\"acklund-Darboux transformation) procedures of explicit solving inverse problem for skew-selfadjoint
Dirac systems developed in \cite{GKS2, KS, FKKS} (see also \cite{SaSaR} and references therein).

In the next section, Preliminaries, we present some basic notions from system theory and formulate several results on Weyl functions.
We also present GBDT  procedure to explicitly solve inverse problem for systems \eqref{1.1}. Namely, we present a procedure to recover skew-selfadjoint Dirac systems (i.e.,
their potentials) from rational Weyl functions (or, more precisely, from minimal realizations of these Weyl functions).
 Section~\ref{Sta} is dedicated to the proof of stability of this procedure. 
 Corollary \ref{Stab1} and Theorem \ref{TmStab2} show stability of the two main steps in solving inverse problem.  Corollary \ref{Stab1} is based on the
 stability of solving the corresponding Riccati equation (see \cite{LaRaT}).
 In addition, new results on so called \cite{GKS2,FKKS}  
pseudo-exponential potentials are formulated
in Lemma \ref{LaSpa} and Corollary \ref{Cyv0}.
Section \ref{discr} is dedicated to the discrete Dirac system. Stability results are presented in Corollary \ref{DStab} and Theorem \ref{TmStab3}.
Uniqueness of the solution of the inverse problem is stated in Theorem \ref{Tmduniq}. Corollary \ref{CyAsC} shows that our sequences $\{C_k\}$
tend to $j$ at infinity. Some proofs are moved into appendix.

As usual, $\BN_0$ stands for the set of non-negative integers, $\BR$ stands for the real axis, $\BC$ stands for the complex plane, $\BC_+$ is the open upper half-plane
$\{z: \,\Im(z)>0\}$, and $\BC_M$ is the open  half-plane $\{ z: \, \Im (z)>M\}$. The notation diag$\{d_1, ...\}$ stands for the diagonal (or block diagonal) matrix 
with the entries $d_1, ...$ on the main diagonal.    By $\s(A)$ we denote   the spectrum of some matrix $A$. (Recall that $\|A\|$ stands for the  $l^2$-induced matrix norm
of $A$.)
We say that the matrix $X$ is positive (nonnegative)
and write $X>0$ ($X\geq 0$) if $X$ is Hermitian, that is,  $X=X^*$, and all the eigenvalues of $X$ are positive (nonnegative).
The notation $I$ stands for the identity operator or matrix and we say that the matrix $X$ is contractive if $ X^*X \leq I$.
Span denotes linear span.
\section{Preliminaries}\label{Prel}
\setcounter{equation}{0}
\subsection{Rational functions}
Recall that  a rational matrix function is called {\it strictly proper} if it tends to zero
at infinity. It is well-known \cite{KFA, LR} that  such an $m_2
\times m_1$ matrix function $\vp$ can be represented in the form
\begin{equation}
\label{app.1} \vp(z)=\clc(z I_n- \cla)^{-1}\clb,
\end{equation}
where $\cla$ is a square matrix of some order $n$, and the matrices $\clb$
and $\clc$ are of sizes $n \times m_1$ and $m_2 \times n$,
respectively. The representation (\ref{app.1})
is called a {\it realization}  of $\vp$, and  the realization
(\ref{app.1}) is said to be {\it minimal} if 
 $n$ is minimal among all possible realizations of $\vp$.
This minimal $n$  is called the {\it McMillan degree} of $\vp$. 
The
realization (\ref{app.1}) of $\vp$ is minimal if and only if
\begin{equation}
\label{app.2} {\mathrm{span}}\bigcup_{k=0}^{n-1}\im
\cla^k\clb=\BC^n,\quad {\mathrm{span}}\bigcup_{k=0}^{n-1}\im  (\cla^*)^k\clc^*=\BC^n, \quad
n= \mathrm{ord}(\cla),
\end{equation}
where Im stands for image and $\mathrm{ord}(\cla)$ stands for the order of $\cla$.
If for a pair of matrices $\{\cla, \, \clb\}$ the first equality in
(\ref{app.2}) holds, then the pair $\{\cla, \, \clb\}$ is called {\it
controllable}. If the second equality in
(\ref{app.2}) is fulfilled, then the pair $\{\clc, \, \cla\}$ is said to be {\it
observable}. 

Differently from the selfadjoint Dirac system case \cite{ALSstab}, where the stability of the solution $X$ of Riccati equation
 $X\clb \clb^*X+\I(\cla^*X-X\cla )+\clc^*\clc =0$
played an important role,
in the case of the skew-selfadjoint Dirac system, we obtain Ricatti equation with minus before $\clb \clb^*$ (see \cite{FKKS} and some references therein):
\begin{align} \label{app.4}&
X\clc^*\clc X+\I(\cla  X-X\cla^*)-\clb \clb^*=0.
\end{align}
From \cite[Proposition 2.2]{GKS2}, which is based on the results from \cite{KFA} (see also \cite[p. 358]{LR} and \cite{FKKS}), we have the statement below.
\begin{Pn} \label{PnRic1} Assume that $\vp(z)$ is a strictly proper rational $m_2 \times m_1$ matrix function and let \eqref{app.1} be its minimal realization.
Then there is a positive solution $X$ $(X>0)$ of the Riccati equation \eqref{app.4}.
\end{Pn}

%%%%%%%%%%%%%%%%%%%%%%%%%%%%%%%%%%%%%%%%%%
%%%%%%%%%%%%%
\subsection{System \eqref{1.1}: Weyl function and inverse problem}\label{Sub2.2}
\begin{Nn}\label{NnFund}
By $Y(x,z)$ we denote the normalized $($by $Y(0,z)=I_m)$ fundamental solution of skew-selfadjoint Dirac system, that is, of system 
\eqref{1.1}, where $j$ and $V$ have the forms \eqref{1.2}.
\end{Nn}
\begin{Dn}  \label{defWeyl1} 
Let Dirac  system \eqref{1.1}, \eqref{1.2} be given and let \eqref{1.2!} hold. Then  an $m_2 \times m_1$ matrix function $\varphi(z)$  such that
\begin{align}&      \label{W1}
\int_0^{\infty}
\begin{bmatrix}
I_{m_1} & \vp(z)^*
\end{bmatrix}
Y(x,z)^*Y(x,z)
\begin{bmatrix}
I_{m_1} \\ \vp(z)
\end{bmatrix}dx< \infty , \quad z \in \BC_M
\end{align}
is called a Weyl function  of the   system \eqref{1.1}, \eqref{1.2}  on $[0, \, \infty)$.
\end{Dn}
Recall that $\BC_M$ is the  half-plane $\{ z: \, \Im (z)>M\}$.
\begin{Rk} \label{Rk2.3} We note that the Weyl function was introduced in \cite{FKRS2012} in an equivalent but different way.
However, Proposition 2.2 and Corollary 2.8 from \cite{FKRS2012} yield the existence and uniqueness of  the function $\vp$ satisfying \eqref{W1}.
This  $\vp(z)$   is holomorphic and contractive in $\BC_M$. 
\end{Rk}
If $\vp$ is rational,  it can be prolonged (from $\BC_M$) on $\BC$  in a natural way.
Each potential $v$ corresponding to a strictly proper rational Weyl function is generated by a fixed value $n\in \BN$
and by a quadruple of matrices, namely, by two $n\times n$ matrices $\a$ and $S_0>0$
and by  $n\times m_k$ matrices $\vt_k$ $(k=1,2)$ such that the matrix identity
\begin{align} \label{2.7}&
\a S_0-S_0\a^*=\I(\vt_1\vt_1^*+\vt_2\vt_2^*) 
\end{align}
holds. Such potentials $v$ have the form
\begin{align} \label{2.5}&
v(x)=2  \vt _{1}^{\, *}\E^{\I x \alpha^{*}} S
(x)^{-1}\E^{\I x
\alpha }\vt_2,
\\ \label{2.6!}& 
S(x)=S_0+ \int_{0}^{x} \Lambda(t) j \Lambda (t)^{*}dt,
\quad
\Lambda (x)= \begin{bmatrix}  \E^{- \I x \alpha } \vt_{1} 
& \E^{\I x
\alpha } \vt_2 \end{bmatrix}.
\end{align}
\begin{Dn}\label{DnAdm}
The quadruples  $\{\a, S_0, \vt_1, \vt_2\}$,
where $S_0>0$ and \eqref{2.7} holds, are called admissible.
\end{Dn}
\begin{Dn} \label{DnPE}\cite{GKS2, FKKS} The potentials $v$, generated $($via equalities \eqref{2.5} and \eqref{2.6!}$)$ by the admissible quadruples  $\{\a, S_0, \vt_1, \vt_2\}$,
are called 
pseudo-exponential.
\end{Dn}
Direct differentiation shows that \eqref{2.7} yields
\begin{align} \label{2.7!}&
\a S(x)-S(x)\a^*=\I \Lam(x)\Lam(x)^*,
\end{align}
that is $\a, \, S(x)$ and $\Lam(x)$ form the so called (see \cite{SaL1, SaL2} and also \cite{SaSaR} and further references therein) $S$-nodes.
\begin{Rk}\label{RkB} According to \cite[Proposition 2.3]{FKKS}, all pseudo-exponential potentials are bounded. Further we show that
pseudo-exponential potentials also tend to zero at infinity.
\end{Rk}
\begin{Tm}\label{TmDp}\cite{FKKS} Let Dirac system \eqref{1.1}, \eqref{1.2} with a pseudo-exponential
potential $v$ be given on $[0,\, \infty)$ and let $v$ be generated by the admissible quadruple $\{\a, S_0, \vt_1, \vt_2\}$.
Then the Weyl function $\vp$ of this system has the form
\begin{align} \label{2.5!}&
\vp(z)=\I\vt_2^*S_0^{-1}(zI_n-\t)^{-1}\vt_1, \quad \t:=\a-\I\vt_1\vt_1^*S_0^{-1}.
\end{align}
\end{Tm}
The following theorem (i.e., \cite[Theorem 2.7]{FKKS}) presents a procedure of explicit solution
of the inverse problem, which is  basic for the next section.

\begin{Tm}\label{TmIpes} Let $\vp(z)$ 
be a strictly proper rational $m_2\times m_1$ matrix function. 
Then $\vp(z)$ is the Weyl function of  the Dirac system \eqref{1.1}, \eqref{1.2} with some pseudo-exponential
potential $v$. This $v$ is uniquely recovered using the following procedure.

Assuming that \eqref{app.1} is a minimal realization of $\vp(z)$ and choosing a positive solution $X>0$
of \eqref{app.4}, we put
\begin{align}
\label{2.6}&
\a=\cla +\I \clb \clb^*X^{-1},  \quad S_0=X, \quad \vt_1=\clb, \quad \vt_2=\I X \clc^*.
\end{align}
The potential $v$ corresponding to the Weyl function $\vp$ is generated $($via \eqref{2.5}, \eqref{2.6!}$)$ by the
quadruple $\{\a, S_0, \vt_1, \vt_2\}$.
\end{Tm}
The matrix identity \eqref{2.7} is immediate from \eqref{app.4} and \eqref{2.6}. Thus, the quadruple constructed in \eqref{2.6}
is admissible. Moreover, each admissible quadruple $\{\a, S_0, \vt_1, \vt_2\}$ satisfies (see \cite[Lemma A.1]{FKKS}) the important relation
\begin{align} \label{sk1}&
\s(\a)\subset (\BC_+\cup \BR).
\end{align}
The quadruples,
which are recovered using \eqref{2.6}, also have an  additional property$:$ controllability of the pair
$\{\a, \vt_1\}$. (This property is immediate from the controllability of the pair
$\{\cla, \clb\}$.) In that case relation \eqref{sk1} may be substituted by a stronger one.
\begin{La}\label{LaSpa} For $\a$ from an admissible quadruple $\{\a, S_0, \vt_1, \vt_2\}$, where \\ $\{\a, \vt_1\}$ is controllable,
we have
\begin{align} \label{sk1!}&
\s(\a)\subset \BC_+ .
\end{align}
\end{La}
\begin{proof}. Putting $\breve \a=S_0^{-1/2}\a S_0^{1/2}$ and $\breve \vt_k=S_0^{-1/2}\vt_k$, we rewrite \eqref{2.7} in the form
\begin{align} \label{2.7+}&
\breve \a -\breve \a^*=\I(\breve \vt_1 \breve \vt_1^*+ \breve \vt_2 \breve \vt_2^*) ,
\end{align}
where $\breve \a$ is linear similar to $\a$. Clearly, the controllability of the pair
$\{\a, \vt_1\}$ yields the controllability of $\{\breve \a, \breve \vt_1\}$.

Assuming that $c\in \BR$ is an eigenvalue of  $\breve \a$, we consider a corresponding eigenvector $g\not=0$
such that $\breve \a g=cg$. Since $c\in \BR$, we obtain $g^*(\breve \a -\breve \a^*) g=0$. Hence, in view of \eqref{2.7+}, we derive 
$$g^*\breve \vt_1=0, \quad g^*\breve \vt_2=0.$$ Therefore, the equalities $\breve \a g=cg$ and \eqref{2.7+} imply that
\begin{align} \label{nk1}&
g^*\breve \a =g^*\breve \a^*=cg^*.
\end{align}
However, the equalities $g^*\breve \vt_1=0$ and $g^*\breve \a =cg^*$ contradict the controllability of the pair $\{\breve \a, \breve \vt_1\}$.
Thus, the relation $\s(\breve \a)\cap \BR=\emptyset$  is proved by negation. Hence, we have $\s( \a)\cap \BR=\emptyset$.
Now, \eqref{sk1!} follows from \eqref{sk1}.
\end{proof}
\begin{Rk}\label{contra}
We note that there are many admissible quadruples generating the same pseudo-exponential potential. 
Furthermore, the matrices $\cla$, $\clb$ and $\clc$ in the minimal realizations \eqref{app.1} of  $\vp$ are unique up to {\rm basis} $($similarity$)$
transformations$:$
\begin{align} \label{2.8}&
\wh{\cla}=\clt^{-1}\cla \clt, \quad \wh{\clc}=\clc \clt, \quad \wh{\clb}=\clt^{-1}\clb,
\end{align}
where $\clt$ are invertible $m\times m$ matrices. 
\end{Rk}

%%%%%%%%%%%%%%%%%%%%%%%%%%%%%%%%%%%%%%%%%%%%%%%%%%%%%
%%%%%%%%%%%%%%%%%%%%%%%%%%%%%%%%%%%%%%%%%%%%%%%%%%%%%%%%%%%%%%%
%%%%%%%%%%%%%%%%%%%%%%%%%%%

\section{Stability}\label{Sta}
\setcounter{equation}{0}
\subsection{Stability of  the recovery of a quadruple}\label{Staq}
First, we consider stability of solving Riccati equation \eqref{app.4}, which appears in Theorem \ref{TmIpes}.
Up to notations, equation \eqref{app.4} coincides with equation (4.1) from \cite{LaRaT}.
\begin{Dn}\cite{LaRaT} \label{Defnong} A nonnegative solution $X$ of \eqref{app.4} is called stably nonnegative
if for every $\ve>0$ there is a $\delta >0$ such that the inequality
\begin{align} \label{nk7}&
\|\cla - \wt{\cla}\|+\| \clb-\wt \clb\|+\|\clc - \wt \clc\|<\delta  
\end{align}
implies that the Riccati equation 
\begin{align} \label{nk8}&
\wt X\wt \clc^* \wt \clc \wt X+\I(\wt \cla \wt  X- \wt X \wt \cla^*)-\wt \clb \wt \clb^*=0
\end{align}
has a nonnegative solution $\wt X$ such that $\|X-\wt X\|<\ve$.
\end{Dn}
Below, we formulate \cite[Theorem 5.4]{LaRaT}, which  describes the stably nonnegative solution of \eqref{app.4}.
\begin{Tm}\label{Tmnong} Assume that the pair $\{\clc, \cla\}$ is observable. Then there is only one
stably nonnegative solution of \eqref{app.4}, being the maximal one.
\end{Tm}
\begin{Rk}\label{Rknong}
The existence of the maximal solution was shown in the earlier papers $($see the discussion at the
beginning of Section 4 in \cite{LaRaT}$)$, and the expression "only one" in Theorem \ref{Tmnong} should
be read as "one and only one". Note that, when $\{\clc, \cla\}$ is observable and also the pair $\{\cla, \clb\}$ is controllable, this maximal solution
is $($in view of  Proposition \ref{PnRic1}$)$ positive. According to  \cite[Theorem 16.3.3]{LR}, this is a unique nonnegative solution as well. 
Thus, $X>0$ considered in Theorem \ref{TmIpes} is unique and  stably positive. 
\end{Rk}
The next corollary follows from  Theorem \ref{Tmnong} and Remark \ref{Rknong}.
\begin{Cy}\label{Stab1}
The recovery $($in Theorem \ref{TmIpes}$)$ of the quadruple $\{\a, S_0, \vt_1,\vt_2\}$ from a triple $\{\cla, \clb, \clc\}$
$($which is
given by a minimal realization \eqref{app.1}$)$ is stable. That is, for every $\ve >0$ there is a $\delta>0$ such that
the inequality $\|\cla -\wt \cla\|+\|\clb -\wt \clb\|+\|\clc -\wt \clc\|<\delta$ yields the inequality
$$\|\a -\wt \a\|+\|S_0 -\wt S_0\|+\|\vt_1 -\wt \vt_1\|+\|\vt_2 -\wt \vt_2\|<\ve,$$
where $\{\wt \a, \wt S_0, \wt \vt_1,\wt \vt_2\}$ is the quadruple corresponding via the procedure from  Theorem \ref{TmIpes}
$($i.e., via the solution $\wt X>0$ of the Riccati equation \eqref{nk8} and via formula \eqref{2.6}$)$ to the minimal realization
$\wt \vp(z)=\wt \clc(z I_n - \wt \cla)^{-1}\wt \clb$.
\end{Cy}

%%%%%%%%%%%%%%%%%%%%%%%%%%%%%%%%%%%%%%%%%%%%%%%%%%%%%%
%%%%%%%%%%%%%%%%%%%%%%%%%%%%%%%%%%%%%%%%%%%%%%%%%%%%%%%%
\subsection{Perturbations of the admissible quadruple}\label{Stap}
Here we will show that small perturbations of the admissible quadruple $\{\a, S_0, \vt_1,\vt_2\}$ result in small
perturbations of the corresponding potential $v$.  For that purpose we will study the matrix function
\begin{align}       & \label{sk2}
R(x)=\E^{-\I x \a}S(x)\E^{\I x \a^*}.
\end{align} 
Expressing $v$ via $R$, we rewrite \eqref{2.5} in the form
\begin{align}       & \label{sk3}
v(x)=2\vt_1^*\E^{2\I x \a^*}R(x)^{-1}\vt_2.
\end{align} 
We note that  only perturbations, which do not change $m_1, m_2$ and $n$, are considered.

It was shown in \cite{ALSstab} that (for the case of selfadjoint Dirac system) a certain matrix function $Q(x)$ monotonically increases to infinity,
and in this subsection we will show that the same holds (in our case) for $R(x)$. However, $R$ differs from $Q$
and the proof is essentially different as well.
\begin{Pn}\label{PnR} Assume that $\{\a, S_0, \vt_1,\vt_2\}$ is an admissible quadruple and that the
pair $\{\a, \vt_1\}$ is controllable. Then $R$ monotonically increases and the minimal eigenvalue of $R(x)$ tends to infinity
$($when $x$ tends to infinity$)$, that is, $R(x)^{-1}$
monotonically decreases and tends to zero.
\end{Pn}
\begin{proof}. Differentiating $R$ and using \eqref{2.6!} and \eqref{2.7!}, we derive
\begin{align}        \label{sk4}
R^{\prime}(x)&=\E^{-\I x \a}\left(S^{\prime}(x)-\I(\a S(x)-S(x)\a^*)\right)\E^{\I x \a^*}
\\ \nn &
=
\E^{-\I x \a}\Lam(x)(j+I_m)\Lam(x)^*\E^{\I x \a^*}=\E^{-2 \I x \a}\vt_1 \vt_1^*\E^{2 \I x \a^*}, \quad R^{\prime}:=\frac{d}{dx}R.
\end{align} 
It easily follows from \eqref{sk2} and \eqref{sk4} that
\begin{align}       & \label{sk5}
R(x)=S_0+\int_0^x\E^{-2 \I t \a}\vt_1 \vt_1^*\E^{2 \I t \a^*}dt.
\end{align} 
Now, it is immediate that $R(x)$ is nondecreasing, and, moreover, $R(x)$ is increasing since $\{\a, \vt_1\}$ is controllable (see \eqref{nk2}).

We prove by negation that the minimal eigenvalue of $R(x)$ tends to infinity (i.e., $R(x)$
tends to infinity).
Indeed, the assumption that  the minimal eigenvalue of $R(x)$ does not tend to infinity implies that there is
a sequence of vectors $g_k$ and values $x_k\in \BR_+$ ($0<k<\infty$) such that $\|g_k\|=1$, $x_k$ tends to infinity (for $k\to \infty$) and 
the sequence $g_k^*R(x_k)g_k$ is bounded. Then there is a partial limit $g\not=0$ of $\{g_k\}$ and for this $g$ we also obtain boundedness:
\begin{align}       & \label{sk6}
\sup_{x\in[0,\infty)}g^*R(x)g<\infty.
\end{align} 
On the other hand, controllability of the pair $\{\a, \vt_1\}$ yields controllability of $\{2 \I \a, \vt_1\vt_1^*\}$. It is well known (see, e.g., \cite{Cop})
that the controllability of  $\{2 \I \a, \vt_1\vt_1^*\}$ is equivalent to the inequality
\begin{align}    &    \label{nk2}
R_0(T):= \int_0^T\E^{-2 \I t \a}\vt_1 \vt_1^*\E^{2 \I t \a^*}dt>0
\end{align} 
for some (and hence for every) $T$. Using \eqref{sk1!} and Jordan normal form of $\a^*$, we can show that for sufficiently large $T>0$ we have
$\| \E^{-2 \I T \a^*}\|\leq 1$ and we fix this $T$. It follows that
\begin{align}    &    \label{nk3}
\| \E^{2 \I T \a^*}f\|\geq \|f\| \quad {\mathrm{for \,\, each}} \quad f\in \BC^n.
\end{align} 
In view of \eqref{nk2} and \eqref{nk3} we obtain
\begin{align}    &    \label{nk4}
f^*R_0(kT)f\geq \ve k f^*f ,\quad 0<k <\infty
\end{align} 
for some $\ve >0$ and for each  $f\in \BC^n$. Since \eqref{nk4}  contradicts \eqref{sk6}, the proposition is proved.
\end{proof}
Matrix identity \eqref{2.7!} together with definition \eqref{sk2} of $R$ and with the second equality in \eqref{2.6!}  imply the identity
\begin{align} \nn&
\a R(x)-R(x)\a^*=\I \E^{-\I x\a}\Lam(x)\Lam(x)^*\E^{\I x\a^*}=\I\big(\E^{-2\I x\a}\vt_1\vt_1^*\E^{2\I x\a^*}+\vt_2\vt_2^*\big).
\end{align}
Multiplying both left and right sides of the identity above by $R^{-1}$ (from the left and from the right), we derive
\begin{align} \label{nk5} &
R(x)^{-1}\a -\a^*R(x)^{-1}=\I R(x)^{-1}\big(\E^{-2\I x\a}\vt_1\vt_1^*\E^{2\I x\a^*}+\vt_2\vt_2^*\big)R(x)^{-1}.
\end{align}
Turning  to the limit in \eqref{nk5},
we see that under conditions of Proposition~\ref{PnR} the following equality holds:
\begin{align} \label{nk6}&
\lim_{x \to \infty} \|\vt_1^*\E^{2\I x\a^*}R(x)^{-1}\|=0.
\end{align}
\begin{Cy}\label{Cyv0}. Each pseudo-exponential potential $v(x)$ tends to zero when $x$ tends to infinity.
\end{Cy}
\begin{proof}. According to Theorems  \ref{TmDp} and \ref{TmIpes}, each pseudo-exponential potential is generated
by some admissible quadruple $\{\a, S_0, \vt_1,\vt_2\}$ such that the
pair $\{\a, \vt_1\}$ is controllable (i.e., the conditions of Proposition \ref{PnR} hold). Now,  our corollary is immediate from 
\eqref{sk3} and \eqref{nk6}.
\end{proof}
We note that only boundedness of $v$ was derived in the previous papers (see Proposition 1.4 \cite{GKS2} and Proposition 2.3 \cite{FKKS}).

The notations corresponding to the quadruples $\{\wt \a, \wt S_0, \wt \vt_1, \wt \vt_2\}$ (in particular, to the perturbed quadruples in the next theorem)
we mark with tilde (e.g., we write $\wt v(x)$,  $\wt R(x)$ and so on).

\begin{Tm} \label{TmStab2}
Let an admissible quadruple $\{\a, S_0, \vt_1,\vt_2\}$, such that the
pair $\{\a, \vt_1\}$ is controllable,  be given. Then, for any $\ve >0$, there is $\delta >0$ such that
each pseudo-exponential potential $\wt v$ generated by an admissible quadruple $\{\wt \a,\wt S_0,\wt \vt_1,\wt \vt_2\}$ satisfying condition
$$\|\a-\wt \a\|+\|S_0-\wt S_0\|+ \|\vt_1-\wt \vt_1\|+\|\vt_2-\wt \vt_2\|<\delta$$
belongs to the $\ve$-neighborhood of $v$ generated by $\{\a, S_0, \vt_1,\vt_2\}$, that is, 
\begin{align}       & \label{3.6}
\sup_{x\in [0,\infty)}\|v(x)-\wt v(x)\|<\ve.
\end{align} 
\end{Tm}
\begin{proof}. Consider  pseudo-exponential potentials $\wt v$ generated by admissible quadruples $\{\wt \a,\wt S_0,\wt \vt_1,\wt \vt_2\}$
belonging to a neighborhood of  $\{\a, S_0, \vt_1,\vt_2\}$.  Recall that the matrix function $R$ corresponding to $\{\wt \a,\wt S_0,\wt \vt_1,\wt \vt_2\}$ is
denoted by $\wt R$. In view of \eqref{sk3}, we have:
\begin{align}     
& \label{3.16}
v(x)=2\vt_1^*\E^{2\I x \a^*}R(x)^{-1}\vt_2, \quad \wt v(x)=2\ \wt \vt_1^*\E^{2\I x \wt \a^*} \wt R(x)^{-1}\wt \vt_2.
\end{align}
Rewriting \eqref{nk5} for $\{\wt \a,\wt S_0,\wt \vt_1,\wt \vt_2\}$, we obtain
\begin{align} \label{nk5!} &
\wt R(x)^{-1}\wt \a -\wt \a^*\wt R(x)^{-1}=\I \wt R(x)^{-1}\big(\E^{-2\I x\wt \a}\wt \vt_1 \wt \vt_1^*\E^{2\I x\wt \a^*}+\wt \vt_2\wt \vt_2^*\big)\wt R(x)^{-1}.
\end{align}
Since conditions of Proposition \ref{PnR} are fulfilled for  $\{\a, S_0, \vt_1,\vt_2\}$, and since $R(x)$ and $\wt R(x)$ are monotonic, we may choose $x_0>0$ and some neighborhood
of  $\{\a, S_0, \vt_1,\vt_2\}$ so that $R(x)$ and $\wt R(x)$ are large enough for $x\geq x_0$. Thus, the left-hand sides of 
\eqref{nk5} and \eqref{nk5!} are small enough. Hence,  the right-hand sides of 
\eqref{nk5} and \eqref{nk5!}  are also small enough. Therefore, taking into account \eqref{3.16}, we see that for any $\ve>0$
there are $x_0>0$ and $\delta_1>0$ such that the next inequality holds in the $\delta_1$-neighborhood 
 of  $\{\a, S_0, \vt_1,\vt_2\}$ (i.e., in the neighborhood  $\|\a-\wt \a\|+\|S_0-\wt S_0\|+ \|\vt_1-\wt \vt_1\|+\|\vt_2-\wt \vt_2\|<\delta_1$):
\begin{align}       & \label{3.18}
\sup_{x\in [x_0,\infty)}\|v(x)-\wt v(x)\|<\ve .
\end{align} 
It easily follows from the definitions of $R$ and $\wt R$ and from \eqref{3.16} that there is some $\delta_2$-neighborhood
of $\{\a, S_0, \vt_1,\vt_2\}$, where we have
\begin{align}       & \label{3.19}
\sup_{x\in [0,x_0)}\|v(x)-\wt v(x)\|<\ve .
\end{align} 
Clearly, inequalities \eqref{3.18} and \eqref{3.19} yield \eqref{3.6}  (for $\delta=\min(\delta_1,\delta_2)$). 
\end{proof}

Corollary \ref{Stab1} and Theorem \ref{TmStab2}
 yield  the stability of the procedure of 
solving inverse problem.

\begin{Cy} The procedure $($given in Theorem \ref{TmIpes}$)$ to uniquely recover the pseudo-exponential potential $v$ of the skew-selfadjoint Dirac system \eqref{1.1} from a minimal
realization of the Weyl function
$($i.e., of some strictly proper rational $m_2\times m_1$  matrix function$)$ is stable.
\end{Cy} 
%%%%%%%%%%%%%%%%%%%%%%%%%%%%%%%%%%%%%%%%%%%%%
%%%%%%%%%%%%%%%%%%%%%%%%%%%%%%%%%%%%%%%%%%%%%%%%%
\section{Discrete Dirac system}\label{discr}
\setcounter{equation}{0}
\subsection{Direct and inverse problems}\label{dDPiIP}
Recall that discrete skew-selfadjoint Dirac system  is given by the formula:
\begin{equation} \label{d1}
y_{k+1}(z)=\left(I_m+ \frac{\I}{ z}
C_k\right)
y_k(z),  \quad C_k=U_k^*jU_k \quad \left( k \in \BN_0 \right).
\end{equation}
\begin{Dn}\label{DnD0} \cite{FKKS} The Weyl function of the discrete system \eqref{d1} is an $m_1 \times m_2$ matrix function $\vp(z)$ in $\BC_M$ $($for some $M>0)$,
which satisfies the inequality
\begin{equation} \label{d1!}
\sum_{k=0}^{\infty}\begin{bmatrix}  \vp(z)^* & I_{m_2}\end{bmatrix}w_k(z)^* w_k(z)
\left[
\begin{array}{c}
 \vp(z) \\ I_{m_2}
\end{array}
\right] < \infty,
\end{equation}
where $w_k(z)$ is the {fundamental solution} of \eqref{d1} normalized by $w_0(z) \equiv I_m$.
\end{Dn}
Similar to the continuous case \eqref{1.1}, the {\it potentials} $\{C_k\}$ of the discrete systems \eqref{d1}
with rational Weyl functions are generated by the admissible quadruples $\{\a, S_0, \vt_1, \vt_2\}$
(see Definition \ref{DnAdm} of the admissible quadruples). More precisely, {\it we additionally require} that 
the pair $\{\a, \vt_1\}$ is controllable, and matrices $C_k$ are determined then by the relations
\begin{align}&\label{d5}
C_k=j+ \Lam_k^* S_k^{-1} \Lam_k - \Lam_{k+1}^* S_{k+1}^{-1} \Lam_{k+1},
\quad k=0,1,2,\ldots;\\
\label{d3}
&\Lam_{k+1}= \Lam_{k}+\I \a^{-1} \Lam_{k} j, \quad  \Lam_{0}=\begin{bmatrix}\vt_1 & \vt_2 \end{bmatrix};
 \\
&
S_{k+1}=S_k+ \a^{-1} S_k (\a^*)^{-1}+ \a^{-1} \Lam_{k} j \Lam_{k}^*
(\a^*)^{-1}.\label{d4}
\end{align}
We note that in the case of the admissible quadruple, where $\{\a, \vt_1\}$ is controllable,
the matrices
 $\a$ and $S_k$ have the following properties (see Lemma \ref{LaSpa} and  \cite[Lemmas 3.2 and  A.1]{FKKS}):
\begin{align}       & \label{dd1}
\s(\a)\in \BC_+, \quad S_k>0 \quad (k\in \BN_0).
\end{align} 
Moreover, in this case, according to \cite[Proposition 3.6]{FKKS}, the matrices $C_k$ given by \eqref{d5}
always admit representation:
\begin{equation} \label{d1+}
 C_k=U_k^*jU_k \quad \left(U_k^*U_k=I_n, \quad k \in \BN_0 \right).
\end{equation}
 The following matrix identities are valid (see \cite[Sect. 3]{FKKS}):
 \begin{align}       & \label{dd1+}
\a S_k-S_k \a^*=\I \Lam_k\Lam_k^* \quad (k\in \BN_0).
\end{align}
According to \cite[Theorem 3.8]{FKKS} the Weyl function $\vp(z)$ of system \eqref{d1}, where $\{C_k\}$ has the form
\eqref{d5}-\eqref{d4}, is given by
\begin{align}       & \label{dd2}
\vp(z)=-\I \vartheta_{1}^{*}S_0^{-1}(z  I_{n} + \g )^{-1}
\vartheta_{2}, \quad \g:=\a-\I \vt_2\vt_2^*S_0^{-1}.
\end{align} 
In particular, $\vp(z)$ is strictly proper rational. Vice versa, given a strictly proper rational $m_1 \times m_2$ matrix function $\vp(z)$
we may recover a system, such that $\vp(z)$ is its Weyl function, using minimal realization \eqref{app.1} (of $\vp(z)$), where
(differently from the continuous case) the matrices $\clb$ and $\clc$ are of sizes  $n\times m_2$ and $m_1 \times n$, respectively.
\begin{Tm}\label{TmDIP} \cite{FKKS} Let $\vp(z)$ 
be a strictly proper rational $m_1\times m_2$ matrix function. 
Then $\vp(z)$ is the Weyl function of  a discrete Dirac system \eqref{d1} with a
potential $\{C_k\}$ generated $($via \eqref{d5}-\eqref{d4}$)$  by some admissible quadruple $\{\a, S_0, \vt_1, \vt_2\}$ 
such that $\{\a, \vt_1\}$ is controllable.
This $\{C_k\}$ is recovered using the following procedure.

Assuming that \eqref{app.1} is a minimal realization of $\vp(z)$ and choosing a positive solution $X>0$
of  the Riccati  equation
\begin{align}
\label{dd3}&
X\clc^*\clc X-\I(\cla  X-X\cla^*)-\clb \clb^*=0.
\end{align}
we put
\begin{align}
\label{dd4}&
\a=- \cla +\I \clb \clb^*X^{-1},  \quad S_0=X,  \quad \vt_1=X \clc^*, \quad \vt_2=\I \clb .
\end{align}
The potential $\{C_k\}$ corresponding to the Weyl function $\vp$ is generated by the
quadruple $\{\a, S_0, \vt_1, \vt_2\}$.
\end{Tm}
\begin{Rk}\label{RkDnong} Theorem \ref{TmDIP} coincides with \cite[Theorem 3.9]{FKKS} after we notice that
the quadruple  $\{\a, S_0, \vt_1, \vt_2\}$  generates the same potential $\{C_k\}$ as the quadruple
$\{X^{-\frac{1}{2}}\a X^{\frac{1}{2}}, I_n, X^{-\frac{1}{2}} \vt_1, X^{-\frac{1}{2}} \vt_2\}$. Moreover, equation \eqref{app.4}
turns into equation \eqref{dd3} if we consider equation \eqref{app.4} corresponding to $-\vp(-z)=\clc(zI_n+\cla)^{-1}\clb$ instead of
\eqref{app.4} corresponding to $\vp(z)$. Hence, since $-\vp(-z)$ is strictly proper rational simultaneously with $\vp(z)$
we can substitute \eqref{dd3} instead of \eqref{app.4} $($as well as $m_2$ instead of $m_1$, and $m_1$ instead of $m_2)$
in Proposition \ref{PnRic1}, Theorem \ref{Tmnong} and Remark \ref{Rknong} and those statements will remain valid.
\end{Rk}
The next corollary follows from  Theorem \ref{TmDIP} and Remark \ref{RkDnong}.
\begin{Cy}\label{DStab}
The recovery $($in Theorem \ref{TmDIP}$)$ of the quadruple $\{\a, S_0, \vt_1,\vt_2\}$ from a triple $\{\cla, \clb, \clc\}$
$($which is
given by a minimal realization \eqref{app.1}$)$ is stable. That is, for every $\ve >0$, there is a $\delta>0$ such that
the inequality $\|\cla -\wt \cla\|+\|\clb -\wt \clb\|+\|\clc -\wt \clc\|<\delta$ yields the inequality
$$\|\a -\wt \a\|+\|S_0 -\wt S_0\|+\|\vt_1 -\wt \vt_1\|+\|\vt_2 -\wt \vt_2\|<\ve,$$
where $\{\wt \a, \wt S_0, \wt \vt_1,\wt \vt_2\}$ is the quadruple corresponding via the procedure from  Theorem \ref{TmDIP}
 to the minimal realization
$\wt \vp(z)=\wt \clc(z I_n - \wt \cla)^{-1}\wt \clb$.
\end{Cy}

\begin{Rk}\label{dcontra} For the quadruple constructed in \eqref{dd4}, the controllability  of  the pair $\{\a,\vt_2\}$ is immediate
$($from the  controllability  of  the pair $\{\cla,\clb\})$, and the controllability  of  the pair $\{\a,\vt_1\}$ follows from the
controllability of the pair $\{X^{-\frac{1}{2}}\a X^{\frac{1}{2}},  X^{-\frac{1}{2}} \vt_1\}$ which is proved in \cite[Appendix A]{FKKS}.
\end{Rk}
According to Remark \ref{dcontra}, the potentials corresponding  to the strictly proper rational Weyl functions
are generated by the quadruples such that the pairs  $\{\a,\vt_1\}$ and $\{\a,\vt_2\}$ are controllable. We introduce
the following definition $($which somewhat differs from the definition in \cite{FKKS}$)$.
\begin{Dn} The quadruple $\{\a, S_0, \vt_1, \vt_2\}$ is called strongly admissible if it is admissible and the pairs
$\{\a,\vt_1\}$ and $\{\a,\vt_2\}$ are controllable. The potentials $\{C_k\}$ generated $($via \eqref{d5}-\eqref{d4}$)$
by the strongly admissible quadruples are called finitely generated and the class of such potentials 
is denoted by the acronym {\rm FG}.
\end{Dn}
The uniqueness of the solution of the inverse problem was not discussed in \cite{FKKS}.
Here, we formulate the uniqueness theorem, the proof of which is given in the appendix.
\begin{Tm}\label{Tmduniq} There is a unique
system \eqref{d1} such that its potential $\{C_k\}$ belongs {\rm FG} and the given 
strictly proper rational $m_1 \times m_2$ matrix function $\vp(z)$ is its Weyl function.
\end{Tm}
This unique system is recovered in Theorem \ref{TmDIP}.
\subsection{Asymptotics of the matrices $C_k$}\label{decr}
Formulas \eqref{d4} and \eqref{dd1+} yield the equality
\begin{align}  \nn
S_{k+1}-&\big(I_n-\I  \a^{-1}\big) S_k\big(I_n+\I ( \a^*)^{-1}\big)
\\ \label{dd6}&
=
\a^{-1}\big(\Lam_k j\Lam_k^*+\Lam_k\Lam_k^*\big)(\a^*)^{-1}.
\end{align} 
For simplicity, we assume further that
\begin{align}
\label{dd5}&
\I \not\in \s(\a).
\end{align}
(Inequality \eqref{dd5} is also essential in the application of the discrete system \eqref{d1} to generalized discrete Heisenberg magnet
model, see \cite{FKKS}.) Then we can introduce the matrices $R_k$:
\begin{align} \label{dd7}
R_k=\big(I_n-\I  \a^{-1}\big)^{-k} S_k\big(I_n+\I ( \a^*)^{-1}\big)^{-k} \qquad (k \geq 0),
\end{align} 
such that (in view of \eqref{dd6}) we have
\begin{align}\nn
 R_{k+1}- R_k=  
(I_{n}-\I  \a^{-1})^{-k-1} & \a^{-1} \Lam_k(I_m+j) \Lam_k^*
\\ \label{dd8} & \,\, 
\times (\a^*)^{-1}(I_{n}+\I ( \a^*)^{-1})^{-k-1} .
\end{align}
Matrices $R_k$ are essential in the study of the asymptotics of $C_k$. Before considering this asymptotics,
we formulate a proposition (which is proved in the appendix) on some interrelations between condition \eqref{dd5}
and finitely generated potentials.
\begin{Pn} \label{Pninotin}
Let an admissible  quadruple  $\{\a, S_0, \vt_1, \vt_2\}$ satisfy the additional relation $0,\I\not\in \s(\a)$.
Then the following statements are valid$:$ \\
$(i)$ We have $S_k>0$ in \eqref{d4}, and so our quadruple generates via \eqref{d5}-\eqref{d4} a well-defined sequence $\{C_k\}$. \\
$(ii)$
There is an admissible quadruple $\{\wt \a, \wt S_0, \wt \vt_1, \wt \vt_2\}$,  which 
satisfies the relation $0,\I\not\in \s(\wt \a)$,
generates the same $\{C_k\}$ as in $(i)$, and has an additional property of controllability of the pair $\{\wt \a,\wt \vt_1\}$. \\
$(iii)$ If we have $C_k\not\equiv j$ for $\{C_k\}$ constructed in $(i)$, there is a strongly admissible
quadruple $\{\wt \a, \wt S_0, \wt \vt_1, \wt \vt_2\}$,  which generates the same $\{C_k\}$ and satisfies the
relation $0,\I\not\in \s(\wt \a)$.
\end{Pn}
%%%%%%%%%%%%%%%%
\begin{Rk} \label{RkI} If the quadruple $\{\wt \a, \wt S_0, \wt \vt_1, \wt \vt_2\}$ satisfies the conditions
given in the statement $(ii)$ or  in the statement $(iii)$ of Proposition \ref{Pninotin}, then
the quadruple $\{\wt S_0^{-\frac{1}{2}}\wt \a \wt S_0^{\frac{1}{2}}, I_n, \wt S_0^{-\frac{1}{2}}\wt \vt_1, \wt S_0^{-\frac{1}{2}}\wt \vt_2\}$
satisfies the same conditions and generates the same potential $\{C_k\}$. Moreover, we have 
\begin{align} \label{dd8+}
\wh S_k=\wt S_0^{-\frac{1}{2}}\wt S_k \wt S_0^{-\frac{1}{2}}, \quad \wh R_k=\wt S_0^{-\frac{1}{2}}\wt R_k \wt S_0^{-\frac{1}{2}},
\end{align} 
where $\wt S_k$ and $\wt R_k$ are generated by the quadruple $\{\wt \a, \wt S_0, \wt \vt_1, \wt \vt_2\}$, and $\wh S_k$ and $\wh R_k$ are generated by the quadruple
$\{\wt S_0^{-\frac{1}{2}}\wt \a \wt S_0^{\frac{1}{2}}, I_n, \wt S_0^{-\frac{1}{2}}\wt \vt_1, \wt S_0^{-\frac{1}{2}}\wt \vt_2\}$.
\end{Rk}
Similar to the case  of the quadruple $\{\wt \a, \wt S_0, \wt \vt_1, \wt \vt_2\}$,  we add  the corresponding accent in  all the notations   
connected with the quadruple $\{\wh \a, \wh S_0, \wh \vt_1, \wh \vt_2\}$    (e.g., we write $\wh S_k$ and $\wh R_k$).

Consider a potential $\{C_k\}$ generated by  some quadruple $\{\a, S_0, \vt_1, \vt_2\}$ satisfying the conditions of Proposition \ref{Pninotin}.
Then, in view of Proposition \ref{Pninotin},  we may assume (without loss of generality) that $\{\a, S_0, \vt_1, \vt_2\}$
is chosen so that it satisfies the conditions  (ii) (on the quadruples $\{\wt \a, \wt S_0, \wt \vt_1, \wt \vt_2\}$) from Proposition \ref{Pninotin},
and we may assume additionally that $S_0=I_n$.

\begin{Pn}\label{PnRk} Let the quadruple $\{\a, S_0, \vt_1, \vt_2\}$ satisfy the conditions (ii) of Proposition \ref{Pninotin}.
Then $R_k$ tends to infinity when $k$ tends to infinity.
\end{Pn}
\begin{proof}.
Let us rewrite \eqref{d3} in an explicit form:
\begin{align} \label{dd9}
\Lam_k=\begin{bmatrix} \big(I_n+\I \a^{-1}\big)^k\vt_1 & \big(I_n-\I \a^{-1}\big)^k\vt_2 \end{bmatrix}.
\end{align} 
Using \eqref{dd9} we rewrite \eqref{dd8} in the form
\begin{equation} \label{dd10} 
 R_{k+1}- R_k=  2(\a-\I I_n)^{-k-1}(\a+\I I_n)^{k}\vt_1\vt_1^*(\a^*-\I I_n)^{k}(\a^*+\I I_n)^{-k-1}.
 \end{equation}
From \eqref{dd10} we derive
\begin{align} \label{dd11} 
 R_{k+n}- R_k= & 2(\a-\I I_n)^{-n-k}(\a+\I I_n)^{k}
\\ \nn & \times
 \left(\sum_{\ell=1}^n (\a-\I I_n)^{n-\ell}(\a+\I I_n)^{\ell-1}\vt_1\vt_1^*(\a^*-\I I_n)^{\ell - 1}(\a^*+\I I_n)^{n-\ell}\right)
 \\ \nn & \times
 (\a^*-\I I_n)^{k}
 (\a^*+\I I_n)^{-n-k}.
 \end{align}
According to Remark \ref{RkI}, we can switch from the quadruple $\{\a, S_0, \vt_1, \vt_2\}$ to the quadruple $\{\wh \a, I_n, \wh \vt_1, \wh \vt_2\}$
and  from $\{R_i\}$ to $\{\wh R_i\}$,
 where
\begin{align} & \label{dd12}
\wh \a=S_0^{-\frac{1}{2}} \a  S_0^{\frac{1}{2}}, \qquad \wh \vt_k= S_0^{-\frac{1}{2}} \vt_k \quad (k=1,2),  \quad
\wh R_k=S_0^{-\frac{1}{2}}R_k S_0^{-\frac{1}{2}};
\\
& \label{dd13}
 \wh R_{k+n}- \wh R_k=S_0^{-\frac{1}{2}}(R_{k+n}- R_k)S_0^{-\frac{1}{2}}.
\end{align}  
The quadruple $\{\wh \a, I_n, \wh \vt_1, \wh \vt_2\}$ has the same properties as $\{\a, S_0, \vt_1, \vt_2\}$ and the operator identity
takes the form
\begin{align} \label{dd14} &
\wh \a - \wh \a^*=\I(\wh \vt _1 \wh \vt_1^*+\wh \vt _2 \wh \vt_2^*).
 \end{align}
 We rewrite \eqref{dd11} for the case of the quadruple $\{\wh \a, I_n, \wh \vt_1, \wh \vt_2\}$:
 \begin{align} \label{dd15} 
 \wh R_{k+n}- \wh R_k= & 2(\wh \a-\I I_n)^{-n-k}(\wh \a+\I I_n)^{k}
\\ \nn & \times
 \left(\sum_{\ell=1}^n (\wh \a-\I I_n)^{n-\ell}(\wh \a+\I I_n)^{\ell-1}\wh \vt_1 \wh \vt_1^*(\wh \a^*-\I I_n)^{\ell - 1}(\wh \a^*+\I I_n)^{n-\ell}\right)
 \\ \nn & \times
 (\wh \a^*-\I I_n)^{k}
 (\wh \a^*+\I I_n)^{-n-k}.
 \end{align} 

Using \eqref{dd14} and \eqref{dd15} we will show that
\begin{align} \label{dd15!} 
 \wh R_{k+n}- \wh R_k  \geq \wh \ve I_n
 \end{align} 
 for some $\wh \ve>0$ which does not depend on $k$.
First, notice that in view of \eqref{dd14} we have
\begin{align} \nn &
(\wh \a-\I I_n)^{-1}(\wh \a+\I I_n)(\wh \a^*-\I I_n)(\wh \a^*+\I I_n)^{-1}-I_n
\\ \label{dd16} & \quad
=2\I(\wh \a-\I I_n)^{-1}(\wh \a^* -\wh \a) (\wh \a^*+\I I_n)^{-1}\geq 0.
 \end{align}
 
 Now, let us show by negation that the sum in the right-hand side of \eqref{dd15} is positive,
 that is, for each vector $f\not= 0$ we have
 \begin{equation} \label{dd17} 
 f^*
 \left(\sum_{\ell=1}^n (\wh \a-\I I_n)^{n-\ell}(\wh \a+\I I_n)^{\ell-1}\wh \vt_1 \wh \vt_1^*(\wh \a^*-\I I_n)^{\ell - 1}(\wh \a^*+\I I_n)^{n-\ell}\right)f\not=0.
 \end{equation} 
 Indeed, if \eqref{dd17} does not hold for some $f\not=0$, we obtain
 \begin{align} &\label{dd18} 
 f^*
\begin{bmatrix} (\wh \a-\I I_n)^{n-1} \wh \vt_1 & (\wh \a-\I I_n)^{n-2}(\wh \a+\I I_n) \wh \vt_1 & \ldots  & (\wh \a+\I I_n)^{n-1}\wh \vt_1 \end{bmatrix}=0.
 \end{align} 
 Recall that the pair $\{\wh \a, \wh \vt_1\}$ is controllable and so
 \begin{align} &\label{dd20} 
{\mathrm{Im}}\, \begin{bmatrix} \wh \a^{n-1} \wh \vt_1 & \wh \a^{n-2} \wh \vt_1 & \ldots  & \wh \vt_1 \end{bmatrix}=\BC^n.
  \end{align}
 It is easy to prove by induction that the polynomials 
 \begin{align} &\label{dd19} 
  (z-\I)^{n-\ell}(z+\I)^{\ell-1} \qquad (1\leq \ell \leq n)
  \end{align} 
 are lineally independent.
 Indeed, if the polynomials 
 $$(z-\I)^{r-\ell}(z+\I)^{\ell-1} \qquad (1\leq \ell \leq r)$$ are lineally independent, then the polynomials
 $$\{(z-\I)(z-\I)^{r-\ell}(z+\I)^{\ell-1}\} \cup (z+\I)^{r+1}$$ 
 are lineally independent as well. Since the polynomials \eqref{dd19} are linearly independent,
 the polynomials $z^{\ell -1}$ can be obtained as linear combinations of the polynomials \eqref{dd19}.
 Therefore (if we switch from $z$ to the matrix $\wh \a$), equality \eqref{dd20} implies that the image of the block matrix in \eqref{dd18} also coincides with $\BC^n$.
 Thus, formula \eqref{dd18} contradicts the condition $f\not=0$  and \eqref{dd17} follows. 
 
 Relations \eqref{dd15}, \eqref{dd16} and \eqref{dd17} yield \eqref{dd15!}. Finally, from  \eqref{dd13} and \eqref{dd15!} we derive
 $R_{k+n}-R_k \geq \ve I_n$ for some $\ve >0$ which does not depend on $k$. Hence, the statement of the proposition is immediate.
 \end{proof}
\begin{Cy} \label{CyAsC}Let the quadruple $\{\a, S_0, \vt_1, \vt_2\}$ satisfy the conditions  of Proposition \ref{Pninotin}. Then, the sequence $C_k$
generated by this quadruple tends to $j$ when $k$ tends to infinity.
\end{Cy}
\begin{proof}.  
Taking into account \eqref{dd1+}, \eqref{dd7} and \eqref{dd9}, we write a matrix identity for $R_k$:
 \begin{equation} \label{dd21} 
 \a R_k -R_k \a^*=\I  \begin{bmatrix} \Psi_k & \vt_2 \end{bmatrix}  \begin{bmatrix} \Psi_k & \vt_2 \end{bmatrix}^*, \quad \Psi_k:= (\a - \I I_n)^{-k} (\a + \I I_n)^{k}\vt_1.
 \end{equation} 
 Moreover, the matrix function $\Lam_k^*S_k^{-1}\Lam_k$ may be written down in the block form in terms of $R_k$ and $\Psi_k$:
   \begin{align} &\label{dd22} 
F_k=\{F_{il k}\}_{i,l=1}^2:=\Lam_k^*S_k^{-1}\Lam_k, \quad F_{11 k}=\Psi_k^*R_k^{-1}\Psi_k, 
\\ &\label{dd23} 
 F_{12 k}=\Psi_k^*R_k^{-1} \vt_2, \quad F_{21 k}=\vt_2^*R_k^{-1}\Psi_k, \quad F_{22 k}=\vt_2^*R_k^{-1}\vt_2.
 \end{align}

 Next, we rewrite \eqref{dd21} in the form
  \begin{align} &\label{dd24} 
R_k^{-1} \a  -\a^*R_k^{-1}=\I R_k^{-1}  \Psi_k \Psi_k^*R_k^{-1}.
 \end{align}
In view of Proposition  \ref{Pninotin}, we assume (without loss of generality) that $\{\a, S_0, \vt_1, \vt_2\}$ satisfies
the conditions (ii) (on the quadruples $\{\wt \a, \wt S_0, \wt \vt_1, \wt \vt_2\}$) from that proposition. Hence, we may apply
Proposition \ref{PnRk}. Proposition~\ref{PnRk} and equality \eqref{dd24} imply that
  \begin{align} &\label{dd25} 
\lim_{k \to \infty} (R_k^{-1}  \Psi_k)=0.
 \end{align}
Partition $C_k$ into the four blocks: $C_k=\{C_{il k}\}_{i,l=1}^2$. Using \eqref{d5}, \eqref{dd23} and \eqref{dd25} we derive
  \begin{align} &\label{dd26} 
\lim_{k \to \infty}C_{12 k} =0, \quad \lim_{k \to \infty}C_{21 k} =0, \quad \lim_{k \to \infty}C_{22 k} =-I_{m_2}.
 \end{align}
 Moreover, it follows from \eqref{d5} that $C_{11 k}=C_{11 k}^*$, and so the matrices $C_{11 k}$ admit representations $C_{11 k}=u_kD_ku_k^*$
 where $u_k$ are unitary matrices and $D_k$ are diagonal. Formulas \eqref{d1+} and \eqref{dd26} yield that
$ \lim_{k \to \infty}C_{11 k}^2=I_{m_1}$, and so $ \lim_{k \to \infty}D_k^2=I_{m_1}$. Hence, relations \eqref{d1+} and \eqref{dd26} yield also that
$D_k>0$ for all sufficiently large $k$. Therefore, the equality $ \lim_{k \to \infty}D_k^2=I_{m_1}$ means that $ \lim_{k \to \infty}D_k=I_{m_1}$.
Thus, we have
  \begin{align} &\label{dd27} 
\lim_{k \to \infty}C_{11 k} =\lim_{k \to \infty}(u_k D_k u_k^*)=I_{m_1}.
 \end{align}
 The statement of the corollary is immediate from \eqref{dd26} and \eqref{dd27}.
\end{proof}
%%%%%%%%%%%%%%%%%%%%%%%%%%%%%%%%%
%%%%%%%%%%%%%%%%%%%%%%%%%%%%%%%%%%
\subsection{Perturbations of the  generating quadruple, \\ discrete case}\label{Stapd}
The following theorem is an analog (for the discrete case) of Theorem \ref{TmStab2}  and has a similar proof.
\begin{Tm} \label{TmStab3}
Let an admissible quadruple $\{\a, S_0, \vt_1,\vt_2\}$, such that the
pair $\{\a, \vt_1\}$ is controllable and $\I \not\in \s(\a)$,  be given. Then, for any $\ve >0$, there is $\delta >0$ such that
each potential $\{\wt C_k\}$ generated by an admissible quadruple $\{\wt \a,\wt S_0,\wt \vt_1,\wt \vt_2\}$ satisfying condition
\begin{align}       & \label{ds1}
\|\a-\wt \a\|+\|S_0-\wt S_0\|+ \|\vt_1-\wt \vt_1\|+\|\vt_2-\wt \vt_2\|<\delta
\end{align}
belongs to the $\ve$-neighborhood of $\{C_k\}$ generated by $\{\a, S_0, \vt_1,\vt_2\}$, that is, 
\begin{align}       & \label{ds2}
\sup_{k\in \BN_0}\| C_k-\wt C_k\|<\ve.
\end{align} 
\end{Tm}
\begin{proof}. Recall that the quadruple $\{\a, S_0, \vt_1,\vt_2\}$ (which is considered in theorem) satisfies \eqref{dd1}, and so
the potential $\{C_k\}$ is well-defined. Without loss of generality we assume that $\delta$ in \eqref{ds1} is sufficiently small, so that the 
admissible quadruples $\{\wt \a,\wt S_0,\wt \vt_1,\wt \vt_2\}$ have the same properties as $\{\a, S_0, \vt_1,\vt_2\}$, namely,
 the pairs $\{\wt \a, \wt \vt_1\}$ are controllable and $\I \not\in \s(\wt \a)$.
We note that in view of \eqref{dd10} (and corresponding formula for $\wt R_k$) the sequences
$\{R_k\}$ and $\{\wt R_k\}$ are nondecreasing. According to the proof of Corollary \ref{CyAsC}, for each $\ve >0$ and each
$\delta$-neighborhood of $\{\a, S_0, \vt_1,\vt_2\}$ we have $\|\wt C_k - j\|<\ve/2$ for sufficiently (depending on $\ve$ and $\delta$)
large values $\wt R_k$. Now, taking into account Proposition \ref{PnRk}, formula \eqref{dd7} and corresponding formula for $\wt R_k$,
we see that we can choose $k_0\in \BN$ and $\delta=\delta_1$ so that $R_k$ and $\wt R_k$ ($k\geq k_0$) are sufficiently large
and so 
\begin{align}       & \label{ds3}
\sup_{k\geq k_0}\| C_k-\wt C_k\|<\ve.
\end{align} 
Moreover,  for each $\ve >0$ and $k_0 \in \BN$ we may choose $\delta_2>0$ so that  
\begin{align}       & \label{ds4}
\sup_{k < k_0}\| C_k-\wt C_k\|<\ve.
\end{align} 
Thus, setting (in \eqref{ds1}) $\delta=\min(\delta_1, \delta_2)$  we derive \eqref{ds2}.
\end{proof}

Corollary \ref{DStab} and Theorem \ref{TmStab3}
 yield  the stability of the procedure of 
solving inverse problem.

\begin{Cy} Consider the procedure of unique recovery of the  potential $\{C_k\}$ of  discrete skew-selfadjoint Dirac system \eqref{d1} from a minimal
realization of the Weyl function $($i.e., of some strictly proper rational $m_1\times m_2$  matrix function$)$, which is given  in Theorem \ref{TmDIP}.
Assume that $\I \not\in \s(\a)$,  where $\a$ is recovered using  \eqref{dd3} and \eqref{dd4}.
Then, this procedure of recovery of the  potential $\{C_k\}$ is stable.
\end{Cy} 
\appendix
\section{Appendix: proofs of Theorem \ref{Tmduniq} \\ and Proposition \ref{Pninotin}}\label{App}
\setcounter{equation}{0}

{\it Proof  of Theorem \ref{Tmduniq}.}
According to Remark \ref{dcontra}, the quadruple $\{ \a, S_0, \vt_1, \vt_2\}$ constructed in Theorem \ref{TmDIP} is strongly admissible,
and so the potential (solution of the inverse problem) $\{C_k\}$,
generated by $\{ \a, S_0, \vt_1, \vt_2\}$,  belongs FG.
Now, assume that another strongly admissible quadruple $\{\wt \a,\wt S_0,\wt \vt_1,\wt \vt_2\}$
generates another potential $\{\wt C_k\}$, such that the Weyl function of the corresponding system is again $\vp$.
Similar to the first phrase in Remark \ref{RkI}, we may assume (without loss of generality) that $S_0=I_n$ and $\wt S_0=I_{\wt n}$.

Since $\vp$ is the Weyl function of  system \eqref{d1} with potential $\{ C_k\}$ and of system \eqref{d1} with potential $\{\wt C_k\}$,
according to \eqref{dd2},  $\vp$  admits two  realizations
\begin{align} \label{ap1}&
\varphi ( z )=-\I  \vartheta_{1}^{*}( z I_{n} +  \g )^{-1}
\vartheta_{2}, \quad \g=\a-\I\vt_2\vt_2^*;
\\ \label{ap2}&
\varphi ( z )=-\I \wt\vartheta_{1}^{*}( z I_{\wt n}+ \wt \g )^{-1}
\wt \vartheta_{2}, \quad \wt \g=\wt \a-\I\wt \vt_2 \wt \vt_2^*.
\end{align}
From the controllability of the pairs $\{\a, \vt_2\}$ and $\{\wt \a,\wt \vt_2\}$, follows
the controllability of the pairs $\{\g, \vt_2\}$ and $\{\wt \g,\wt \vt_2\}$, respectively.
Taking into account 
the controllability of the pairs $\{\a, \vt_1\}$ and $\{\wt \a,\wt \vt_1\}$ and the equalities
$$\a=\a^*+\I(\vt_1\vt_1^*+\vt_2\vt_2^*)=\g^*+\I \vt_1\vt_1^*, \quad 
\wt \a=\wt \a^*+\I(\wt \vt_1\wt \vt_1^*+\wt \vt_2\wt \vt_2^*)=\wt \g^*+\I \wt\vt_1\wt\vt_1^*,$$
we derive the controllability of the pairs $\{\g^*, \vt_1\}$ and $\{\wt \g^*,\wt \vt_1\}$,
that is, the observability of the pairs $\{ \vt_1^*,  \g\}$ and $\{\wt \vt_1^*, \wt \g\}$.
The controllability of the pairs $\{\g, \vt_2\}$ and $\{\wt \g,\wt \vt_2\}$ and
the observability of the pairs $\{ \vt_1^*,  \g\}$ and $\{\wt \vt_1^*, \wt \g\}$ yield
the minimality of  both realizations \eqref{ap1} and \eqref{ap2}. (In fact, the minimality
of \eqref{ap1} follows from the proof of Theorem \ref{TmDIP}.)

Our next arguments coincide with the final arguments in the proof of Theorem 0.5 \cite{KS}.
Since the realizations \eqref{ap1} and \eqref{ap2} are minimal, we have $n=\wt n$ 
and there exists a nonsingular matrix $\cls$ such that
\begin{align} \label{ap3}&
\wt \g=\cls\g\cls^{-1}, \quad \wt\vt_2=\cls\vt_2, \quad \wt \vt_1^*=\vt_1^*\cls^{-1}
\end{align}
(see Remark \ref{contra}). Identities \eqref{2.7} for our quadruples may be rewritten:
\begin{align} \label{ap4}&
\g -\g^*=\I(\vt_1\vt_1^*-\vt_2\vt_2^*), \quad \wt \g -\wt \g^*=\I(\wt \vt_1 \wt \vt_1^*-\wt \vt_2 \wt \vt_2^*).
\end{align}
Substituting \eqref{ap3} into the second equality in \eqref{ap4} and multiplying the result
by $\cls^{-1}$ from the left and by $(\cls^*)^{-1}$ from the right, we obtain
\begin{align} \label{ap5}&
\g Z  -Z \g^*=\I(Z\vt_1\vt_1^*Z- \vt_2\vt_2^*), \quad Z:=(\cls^*\cls)^{-1}.
\end{align}
Due to the controllability of the pairs $\{\g, \, \vt_2\}$ and   $\{\g^*, \, \vt_1\}$ and to the second part of Remark \ref{Rknong},
there is a unique nonnegative solution of \eqref{ap5}. 
Therefore, comparing \eqref{ap5} and the
first equality in \eqref{ap4}, we see that $Z=I_n$, that is, $\cls$ is unitary. 
Taking into account that $\cls$ is unitary, we easily derive from \eqref{d5}--\eqref{d4} and from
the state-similarity transformation \eqref{ap3} the equality $C_k\equiv \wt C_k$.
The uniqueness of the solution of the inverse problem follows.  $\blacksquare$

\vspace{0.3em}

{\it Proof  of  Proposition \ref{Pninotin}.}\\
Step 1. 
In view of $0,\I \not\in \s(\a)$, the statement of (i) (i.e., the inequality $S_r>0$ for $r\geq 0$) follows by induction
from $S_0>0$ and from \eqref{dd6}.
Moreover, similar to the proof of Theorem \ref{Tmduniq} 
we may restrict further proofs to the case $S_0=I_n$ and $\wt S_0 =I_{\wt n}$.

Let the admissible quadruple $\{ \a, I_n, \vt_1, \vt_2\}$, such that  $0,\I \not\in \s(\a)$, generate $\{C_k\}$ but assume that
 the pair $\{\a, \, \vt_1\}$ is not controllable. Put
\begin{align} &\label{d48}
 L_0: =  \spa
\bigcup_{k=0}^{n-1} \im  \big(\a^k  \vt_1), \quad
\wt n:= \dim  L_0 , \quad \wt L_0:=  \im \begin{bmatrix}0\\ I_{\wt n} \end{bmatrix}
\in \BC^n,
\end{align} 
and consider, first, the case $\wt n>0$. For that case we choose a unitary matrix $q$ that
maps $ L_0$ onto  $\wt L_0$ and consider matrices
\begin{align} &\label{ad1}
\wh \a:=q  \a q^*, \quad \wh \vt_1:= q  \vt_1, \quad \wh \vt_2:= q  \vt_2.
\end{align} 
Let us show that these  matrices have the following block structure:
\begin{equation}\label{da9}
\wh \a=
\left[ \begin{array}{lr}\wh  \a_{11} & 0 \\
\wh  \a_{21} & \wt \a  \end{array} \right], \quad \wh \vt_1= \left[
\begin{array}{c}0 \\ \wt  \vt_{1}
\end{array} \right], \quad
\wh \vt_2=
\left[ \begin{array}{c} \vk \\ \wt \vt_{2}
 \end{array} \right].
\end{equation}
Indeed, since
$L_0$ is an invariant subspace of $ \a$, we see that $\wt L_0$ is an
invariant subspace of $ \wh \a$, and thus $\wh \a$ has the block
triangular form given in (\ref{da9}). Next,  notice that the
inclusion $\im( \vt_1) \subseteq  L_0$ yields  $\im(\wh \vt_1)
\subseteq \wt  L_0$, that is,
$\wh \vt_1$ also has the block form given in
(\ref{da9}). 

Taking into account that $q$ is unitary, that
 the quadruple $\{ \a, I_n, \vt_1, \vt_2\}$ is admissible and that $0,\I \not\in \s(\a)$,
 we see that the quadruple $\{\wh \a, I_n, \wh \vt_1, \wh \vt_2\}$ has the same properties.
 Moreover, in view of \eqref{da9},  the quadruple $\{\wt \a, I_{\wt n}, \wt \vt_1, \wt \vt_2\}$
 has the same properties as well.

 The quadruple $\{ \a, I_n, \vt_1, \vt_2\}$
determines (via formulas \eqref{d3} and \eqref{d4}) the matrices
$\Lam_k$ and $S_k$. Finally, formula \eqref{d5} determines $C_k$.
As before, we use the accents ``hat" and ``tilde" 
  in the notations
of the matrices (e.g., of $\Lam_k$, $S_k$, $R_k$ and $C_k$) determined by the 
quadruples $\{\wh \a, I_n, \wh \vt_1, \wh \vt_2\}$ and $\{\wt \a, I_{\wt n}, \wt \vt_1, \wt \vt_2\}$, respectively.
For instance, we write $\wh C_k$ and $\wt C_k$.   
It is immediate that $\wh C_k \equiv C_k$.
In order to prove that $\wt C_k \equiv \wh C_k$ we make some preparations.

Rewriting \eqref{dd7} for the cases of $\wh R_{k+1}-\wh R_k$ and $\wt R_{k+1}-\wt R_k$ and using \eqref{da9} we derive
 \begin{align}&  
\wh R_{k+1}-\wh R_k =  
\begin{bmatrix} 0 & 0 \\ 0 & \wt R_{k+1}-\wt R_k 
\end{bmatrix}.
 \nonumber
\end{align}
Hence, taking into account \eqref{dd7} (for $k=0$), we have
 \begin{align}&  \label{da17}
\wh R_k = {\mathrm{diag}} \, \{ I_{n-\wt n}, \, \wt R_k
 \}.
\end{align}
Rewriting formulas \eqref{dd22}, \eqref{dd23} for the cases of $\wh \Lam_k^* \wh S_k^{-1} \wh \Lam_k$ and $\wt \Lam_k^* \wt S_k^{-1} \wt \Lam_k$
(instead of $ \Lam_k^*  S_k^{-1}  \Lam_k$), according to \eqref{da9} and \eqref{da17} we obtain
\begin{equation}\label{da18}
\wh \Lam(k)^*\wh S_k^{-1}\wh \Lam(k)- \wt \Lam(k)^* \wt S_k^{-1} \wt \Lam(k)= \left[
\begin{array}{lr}
0 & 0 \\ 0 & \varkappa^* \varkappa
\end{array}
\right].
\end{equation}

Finally, relations (\ref{d5}) and (\ref{da18}) yield the equality $\wt C_k\equiv \wh C_k$. Recalling that
$\wh C_k\equiv  C_k$, we see that the quadruple $\{\wt \a, I_{\wt n}, \wt \vt_1, \wt \vt_2\}$
generates $\{C_k\}$. Moreover, this quadruple has all the properties of $\{ \a, I_n, \vt_1, \vt_2\}$ and the pair $\{\wt \a, \wt \vt_1 \}$ is controllable.

In order to complete the proof of (ii), it remains to consider the case $\wt n =0$. Clearly, $\wt n =0$ implies that $\vt_1=0$. Therefore,
formulas \eqref{dd7}, \eqref{dd8} and \eqref{dd9} yield $R_k\equiv I_n$. Since $\vt_1=0$ and $R_k\equiv I_n$,  it follows from \eqref{dd21}-\eqref{dd23}
that 
\begin{align}&  \label{da19}
\Lam_k^*S_k^{-1} \Lam_k \equiv \left[
\begin{array}{lr}
0 & 0 \\ 0 & \vt_2^*\vt_2
\end{array}
\right], \quad \Lam_k^*S_k^{-1} \Lam_k-\Lam_{k+1}^*S_{k+1}^{-1} \Lam_{k+1} \equiv 0.
 \end{align} 
From \eqref{d5} and \eqref{da19}, it is immediate that $C_k\equiv j$. Hence, we can choose
any admissible quadruple, such that  $\I \not\in \s(\a)$, the pair $\{\a, \vt_1\}$ is controllable   and $\vt_2=0$. 
Let us show that such quadruple  will generate our
$\{C_k\}$. 
Indeed, recall that \eqref{dd1} follows from the controllability of $\{\a, \vt_1\}$. So, $0, - \I \not\in \s(\a)$. Similar to the proof of \eqref{dd8}, we can show that the equality
\begin{align}\nn
 Q_{k+1}- Q_k =& \big(I_{n}+\I  \a^{-1}\big)^{-k-1}\a^{-1} 
\\ \label{da20}   &\,\, \times
\Lam(k)(j-I_m) \Lam(k)^* (
\a^*)^{-1}\big(I_{n}-\I ( \a^*)^{-1}\big)^{-k-1} 
\end{align} 
is valid for
 \begin{align} \label{da14}
Q_k:=\big(I_n+\I  \a^{-1}\big)^{-k} S_k\big(I_n-\I ( \a^*)^{-1}\big)^{-k}.
\end{align} 
Then, we express $\Lam_k^*S_k^{-1} \Lam_k$ in terms of $Q_k$ instead of $S_k$, and the proof
that an admissible quadruple with $\vt_2=0$ generates $C_k\equiv j$ is similar to the proof (above) that 
an admissible quadruple with $\vt_1=0$ generates $C_k\equiv j$. 

The statement (ii) is proved.

 Step 2.   In order to prove $(iii)$, we consider again some potential  $\{C_k\}$ $(C_k \not\equiv j)$ generated by the quadruple $\{ \a, I_n, \vt_1, \vt_2\}$
 satisfying conditions of Proposition \ref{Pninotin}. Without loss of generality, we assume that the quadruple $\{ \a, I_n, \vt_1, \vt_2\}$ is chosen so
 that $n$ there is the least value of $n$ for the quadruples generating $\{C_k\}$ and satisfying conditions of Proposition \ref{Pninotin}.
 In view of the statement (ii) of this proposition and the corresponding constructions in Step 1,
 the pair $\{\a, \vt_1\}$ is controllable. Assuming that   the pair $\{\a, \, \vt_2\}$ is not controllable, we construct below
a quadruple $\{\wt \a, I_{\wt n}, \wt \vt_1, \wt \vt_2\}$ generating $\{C_k\}$ and satisfying the conditions of Proposition \ref{Pninotin}, where  $\wt n < n$.
Thus, we come to a contradiction, which implies the controllability of $\{\a, \, \vt_2\}$.

First, put
\begin{align} &\label{d48!}
 L_0: =  \spa
\bigcup_{k=0}^{n-1} \im  \big(\a^k  \vt_2), \quad
\wt n:= \dim  L_0 , \quad \wt L_0:=  \im \begin{bmatrix}0\\ I_{\wt n} 
 \end{bmatrix} \in \BC^n,
\end{align}  
and let $\wt n <n$. Since $C_k \not\equiv j$, the case $\wt n =0$ is excluded, because
 in that case we have $\vt_2=0$, which yields $C_k\equiv j$ (see the arguments at the end of Step 1).
 The following considerations are also similar to the considerations in Step 1, although we deal with
 the matrices $Q_k$ instead of the matrices $R_k$.
 Namely, we choose a unitary matrix $q$ which
maps $ L_0$ onto  $\wt L_0$. We introduce matrices $\wh \a, \wh \vt_1, \wh \vt_2$ via formula \eqref{ad1}.
The structure of these matrices may be proved in the same way as the formula \eqref{da9} and we obtain:
\begin{equation}\label{2da9}
\wh \a=
\left[ \begin{array}{lr}\wh  \a_{11} & 0 \\
\wh  \a_{21} & \wt \a  \end{array} \right], \quad \wh \vt_1= \left[
\begin{array}{c}\vk \\ \wt  \vt_{1}
\end{array} \right], \quad
\wh \vt_2=
\left[ \begin{array}{c} 0 \\ \wt \vt_{2}
 \end{array} \right].
\end{equation}
We note that again  $\wh C_k \equiv C_k$ and the quadruple $\{\wt \a, I_{\wt n}, \wt \vt_1, \wt \vt_2\}$
satisfies the conditions of Proposition \ref{Pninotin}.

It remains to show that $\wt C_k \equiv \wh C_k$. Rewriting \eqref{da20} for the cases of the
quadruples $\{\wh \a,I_n, \wh \vt_1, \wh \vt_2\}$ and $\{\wt \a, I_{\wt n}, \wt \vt_1, \wt \vt_2\}$
and taking into account \eqref{2da9},
we obtain
 \begin{align} &\label{da22}
\wh Q_{k+1}-\wh Q_k=\begin{bmatrix}0 & 0 \\ 0 & \wt Q_{k+1}-\wt Q_k\end{bmatrix}.
\end{align}  
Thus, the matrices $\wh Q_k$ are block diagonal and have the form 
 \begin{align} &\label{da23}
\wh Q_k=\diag\{I_{n-\wt n}, \, \wt Q_k\}.
\end{align}  
Substituting
 \begin{align} \label{2da14}
\wh S_k=\big(I_n+\I \wh \a^{-1}\big)^{k}\wh Q_k\big(I_n-\I (\wh \a^*)^{-1}\big)^{k}
\end{align} 
into $\wh \Lam_k^*\wh S_k^{-1} \wh\Lam_k$ and using \eqref{dd9} (rewritten for $\wh\Lam_k$),
  \eqref{2da9} and \eqref{da23}, we derive
the identity
\begin{equation}\label{da24}
\wh \Lam_k^*\wh S_k^{-1}\wh \Lam_k- \wt \Lam_k^* \wt S_k^{-1} \wt \Lam_k= \left[
\begin{array}{lr}
\varkappa^* \varkappa & 0 \\ 0 & 0
\end{array}
\right].
\end{equation}
Equalities \eqref{d5} and \eqref{da24} show that $\wt C_k\equiv \wh C_k \,(=C_k)$.
Thus, the sequence $\{C_k\}$ is generated by the  quadruple $\{\wt \a, I_{\wt n}, \wt \vt_1, \wt \vt_2\}$ ($\wt n <n$)   and 
we arrive at a contradiction. The statement (iii) is proved.
$\blacksquare$
%%%%%%%%%%%%%%%%

%%%%%%%%%%%%%%%%%%%%%%%%%%%%%%%%%%%%%%%%%%%%%%%%%%%%%%%%%%%%%%%%%%
%%%%%%%%%%%%%%%%%%%%%%%%%%%%%%%%%%%%%%%%%%%%%%%%%%%%%%%%%%
%%%%%%%%%%%%%%%%%%%%%%%%%%%%%%%%%%

\bigskip 
\noindent{\bf Acknowledgments.}
 {The research of A.L. Sakhnovich   was supported by the
Austrian Science Fund (FWF) under Grant  No. P24301.}
%%%%%%%%%%%%%%%%%%%%%%%%%%%%%%%%%%%%%%%%%%%%%%
%%%%%%%%%%%%%%%%%%%%%%%%%%%%%%%%%%%%%%%%%%%%%%%
\newpage
%%%%%%%%%%%%%%%%%%%%%%%%%%%%%%%%%%%%%%%%%%%
%%%%%%%%%%%%%%%%%%%%%%%%%%%%%%%%%%%%%%%%%%%%

\begin{flushright}
B. Fritzsche,\\
Fakult\"at f\"ur Mathematik und Informatik,  Universit\"at Leipzig, \\
 Augustusplatz 10,  D-04009 Leipzig, Germany, \\
 e-mail: {\tt Bernd.Fritzsche@math.uni-leipzig.de}
 
\vspace{0.5em}

B. Kirstein, \\ 
Fakult\"at f\"ur Mathematik und Informatik,  Universit\"at Leipzig, \\
 Augustusplatz 10,  D-04009 Leipzig, Germany, \\
 e-mail: {\tt Bernd.Kirstein@math.uni-leipzig.de}
 
\vspace{0.5em}

I.Ya. Roitberg, \\
 e-mail: {\tt  	innaroitberg@gmail.com}

\vspace{2em} 

A.L. Sakhnovich,\\
Institute for Analysis and Scientific Computing,\\
Vienna
University
of
Technology, \\
Wiedner Haupstr. 8-10 / 101, 1040 Wien,
Austria, \\
e-mail: {\tt oleksandr.sakhnovych@tuwien.ac.at}

\end{flushright}

%%%%%%%%%%%%%%%%%%%%%%%%

\end{document}